\documentstyle%{article}
[12pt]{article}
\newtheorem{df}{ \sc Definition}[section]
\newtheorem{ex}[df]{ \sc Example}
\newtheorem{pr}[df]{ \sc Proposition}
\newtheorem{th}[df]{ \sc Theorem}

\newtheorem{re}[df]{ \sc Remark}
\newtheorem{lem}[df]{ \sc Lemma}
\def\i{{\bf i}^*}
\def\j{{\bf j}^*}
\def\k{{\bf k}^*}
\def\ii{{\bf i}}
\def\o{\otimes}
\def\M{M^\x}
\def\x{{^{_\bullet}}}
\def\g{{\bf g}}
\def\L{\Lambda^\x\hskip -4pt\g^*}
\def\Lg{\Lambda^\x\hskip -4pt\g}
\def\Sg{S^\x\hskip -4pt\g^*}
\def\Ll{{\cal L}}
\def\dl{d_\Lambda}
\def\mono{\lhook\joinrel\relbar\joinrel\rightarrow}
\def\epi{\mathrel -\joinrel\rightarrow\hskip -10pt\rightarrow}
\begin{document}

\author{Tomasz Maszczyk and Andrzej Weber\footnote{Supported by KBN
2P03A 00218 grant.}\\
Institute of Mathematics, Uniwersytet Warszawski \\
 ul. Banacha 2, 02--097 Warszawa, Poland\\
e-mail:maszczyk@mimuw.edu.pl, aweber@mimuw.edu.pl }
\title{ Koszul Duality for modules over Lie algebras }
\maketitle \noindent\hfil {Math. Sub. Class.: Primary 17B55,
17B20; Secondary 55N91}\vskip 10pt

Let $\g$ be a reductive Lie algebra over a field of characteristic
zero. Suppose $\g$ acts on a complex of vector spaces $\M$ by
$i_\lambda$ and $\Ll_\lambda$, which satisfy the same identities that
contraction and Lie derivative do for differential forms. Out of
this data one defines the cohomology of the invariants and the equivariant
cohomology of $\M$. We establish Koszul duality between each
other.

\section {Introduction}
Let $G$ be a compact Lie group.  Set $\Lambda_\x=H_*(G)$ and
$S^\x=H^*(BG)$. The coefficients are in $\bf R$ or $\bf C$.
Suppose $G$ acts on a reasonable space $X$.  In the paper \cite{GKM}
Goresky, Kottwitz and MacPherson established a duality between the
ordinary cohomology which is a module over $\Lambda_\x$ and
equivariant cohomology which is a module over $S^\x$. This duality is
on the level of chains, not on the level of cohomology.
Koszul duality says that there is an equi\-va\-lence of
derived categories of $\Lambda_\x$--modules and $S^\x$--modules.
One can lift the structure of an $S^\x$--module on
$H_G^*(X)$ and the structure of a $\Lambda_\x$--module on $H^*(X)$
to the level of chains in such a way that the obtained complexes
correspond to each other under Koszul duality. Equivariant
coefficients in the sense of \cite{BL} are also allowed. Later,
Allday and Puppe (\cite{AP}) gave an explanation for this duality
based on the minimal Brown-Hirsh model of the Borel construction. One
should remark that Koszul duality is a reflection of a more general
duality: the one described by Husemoller, Moore and Stasheff in
\cite{HMS}.

Our goal is to show that this duality phenomenon is a purely algebraic
affair. We will construct it without appealing to topology. We
consider a reductive Lie algebra $\g$ and a complex of vector
spaces $\M$ on which $\g$ acts via two kinds of actions:
$i_\lambda$ and $\Ll_\lambda$. These actions satisfy the same
identities as contraction and Lie derivative do in the case of the
action on the differential forms of a $G$--manifold. Such
differential $\g$--modules
were already described by Cartan in \cite{Ca}; see also \cite{AM},
\cite{GS}. We do not assume that $\M$ is finite dimensional nor
semisimple. We also wish to correct a small inaccuracy in the
proof of Lemma 17.6, \cite{GKM}. The distinguished transgression
plays a crucial role in our construction. This is a canonical
identification between the space of primitive elements of
$(\L)^\g\simeq H^*(\g)$ with certain generators of $(S\g^*)^\g$.

The results of the present paper were obtained with the help and
assistance of Marcin Cha\l upnik. It is a part of joint work. Our
aim is to describe Koszul duality in a much wider context.
We thank the referee for careful corrections.

%{\def\contentsname{\normalsize\bf Contents:}\tableofcontents}

\section{Category }
 Let $\g$ be a reductive Lie algebra over a field $k$ of
characteristic 0. We consider differential graded vector spaces
$\M $ over $k$ equipped with linear operations $i_\lambda:\M
\rightarrow M^{\x-1}$ of degree $-1$ for each $\lambda\in\g$. We
define $$\Ll_\lambda=di_\lambda+i_\lambda d:\M \rightarrow \M
\,.$$ We assume that $i_\lambda$ is linear with respect to
$\lambda$ and for each $\lambda, \mu \in \g$ the following
identities are satisfied: $$i_\lambda i_\mu = -i_\mu i_\lambda,$$
$$[\Ll_\mu,i_\lambda]=i_{[\mu,\lambda]}\,.$$ Then $\Ll$ is a
representation of $\g$ in $\M $. The category of such objects with
obvious morphisms will be denoted by $K(\g)$.

\begin{ex} \rm Let $G$ be a group with Lie algebra $\g$ and let $X$ be
a $G$-manifold. Then the space of differential forms
$\Omega^\x(X)$ equipped with the contractions
with fundamental vector fields is an example of an object from
$K(\g)$.\end{ex}

\begin{ex} \rm Another example of an object of $K(\g)$ is $\L$, the
exterior power of the dual of $\g$. The generators of $\L$
are given the gradation 1. This is a differential graded algebra
with a differential $\dl$ induced by the Lie bracket. The operations
$i_\lambda$ are the contractions with $\g$.\end{ex}

\begin{ex} \rm For a representation $V$ of $\g$ define the invariant
subspace $$V^\g=\{v\in V:\;\forall \lambda\in \g \;\; \Ll_\lambda
v =0\}.$$ Then $V^\g$ with trivial differential and $i_\lambda$'s
is an object of $K(\g)$. In particular we take $V=\Sg$,  the
symmetric power of the dual of $\g$. The generators of $\Sg$
are given the gradation 2.\end{ex}

\begin{ex} \rm Suppose there are given two objects $\M $ and $N^\x$ of
$K(\g)$. Then $\M \o N^\x$ with operations $i_\lambda$ defined by the
Leibniz formula is again in $K(\g)$. \end{ex}

Note that the objects of $K(\g)$ are the same as differential
graded modules over a dg--Lie algebra $C\g$ (the cone over $\g$)
with $$C\g^0=C\g^{-1}=\g\qquad  C\g^{\neq -1,0}=0\,.$$ The
elements of $C\g^0$ are denoted by $\Ll_\lambda$ and the elements
of $C\g^{-1}$ are denoted by $i_\lambda$. They satisfy the following
identities:
$$di_\lambda=\Ll_\lambda\,,
\qquad[\Ll_\lambda,\Ll_\mu]=\Ll_{[\lambda,\mu]}\,,
\qquad[\Ll_\lambda,i_\mu]=i_{[\lambda,\mu]}\,,
\qquad[i_\lambda,i_\mu]=0\,.$$

The enveloping dg--algebra of $C\g$ is the Chevaley--Eilenberg
complex $V(\g)$ (\cite{Wei} p. 238) which is a free
$U(\g)$--resolution of $k$. Thus the objects of $K(\g )$ are just
the dg--modules over $V(\g)$.

\begin{ex} \rm Suppose that a Lie algebra $\g$ acts on  a
 graded commutative $k$-algebra $A$ by derivations. Let $\Omega_k A$ be the
 algebra of forms; it is generated by symbols $a$ and $da$ with
 $a\in A$. The cone of $\g$ acts
 by derivations. The action on generators is given by
 $$ i_\lambda a=0,\quad \Ll_\lambda a=\lambda a,
 \quad i_\lambda da=\lambda a,\quad \Ll_\lambda da=d\lambda
 a\,.$$\end{ex}
 Then $\Omega_k A$ is in $K(\g)$.

 Another point of view (as in \cite{GS}) is that the objects of
$K(\g)$ are the representations of the super Lie algebra
$\widehat\g=C\g\oplus k[-1]$ (where $k[-1]$ is generated by $d$)
with relations: $$[d,d]=0,\qquad [d,x]=d x\qquad{\rm for}\; x\in
C\g\,.$$

\section{The twist}

 We will describe a transformation, which plays the role of the
canonical map $$\begin{array}{ccccc}\Phi:& G\times X&\rightarrow&
G\times X\\ & (g,x)&\mapsto &(g,gx)\\
\end{array}$$
for topological $G$-spaces.

  We define a linear map $\ii:\L\o\M\rightarrow\L\o\M$ by the
 formula:
$$\ii(\xi\o m)=\sum_k\xi\wedge\lambda^k\o i_{\lambda_k}m\,,$$
where $\{\lambda_k\}$ is a basis of $\g$ and $\{\lambda^k\}$ is
the dual basis. It commutes with $\Ll_\lambda$.
The operation $\ii$ is nilpotent. Define an automorphism
 ${\bf T}:\L\o\M\rightarrow\L\o\M$ (which is not in $K(\g)$):
$${\bf T}=\exp(-\ii)=\sum_{n=1}^{\infty}\frac{(-1)^n}{n!}\ii^n\,,$$
$${\bf T}(\xi\o m)=
\sum_{I=\{i_1<\dots<i_n\}}(-1)^{\frac{n(n+1)}2}
\xi\wedge\lambda^I\o i_{\lambda_I}m\,.$$
 It satisfies (\cite{GHV}, Prop.~V, p.286, see also \cite{AM})
$$i_\mu({\bf T}(\xi\o m))={\bf T}((i_\mu\xi)\o m)\,,$$ $$d\,{\bf T}
(\xi\o m) = {\bf T}\left( d(\xi\o m)+
\sum_k\lambda^k\wedge\xi\o\Ll_{\lambda_k}m\right)\,.$$

\begin{re} \rm Note that for the self-map $\Phi$ of $G\times X$ we have
$$d\Phi^*\omega(1,x)=\left(d\omega+\sum_k
p^*\lambda^k\wedge\Ll_{\lambda_k}\omega\right)(1,x)\,,$$ where
$p:G\times X\rightarrow X$ is the projection and  $x\in X$.\end{re}

The twist on the level of the Weil algebra has already been used by
Cartan \cite{Ca} and later by Mathai and Quillen \cite{MQ}.
From another point of view, for a d.g.vector space to be an object of $K(\g)$
is equivalent to having such a twist which satisfies certain
axioms. We will not state them here. We just remark that 
understanding of this twist allows to develop a theory of actions of
${\cal L}_\infty$-algebras.

\section{Weil algebra }

 Following \cite{GHV}, Chapter VI, p.223 we define the Weil
algebra $$W(\g)=\Lambda^\x(C\g )^*=\Sg \o\L\,.$$ The generators of
$\Sg$ are given the gradation 2 whereas the generators of $\L$ are
given the gradation 1. The differential in $W(\g)$ is the sum of
three operations: $$ d_W(a\o b)  =a\o\dl b+
  \sum_k \lambda^k a\o i_{\lambda_k}b
 + \sum_k ad^*_{\lambda_k} a\o \lambda^k\wedge b\,,$$
where $\{\lambda_k\}$ is a basis of $\g$.

The differential $d_W$ satisfies Maurer--Cartan formula
$$d_W(1\o \xi)-1\o \dl\xi = \xi\o 1\,,$$
for $\xi\in\g^*$.

 The operations $i_\lambda$ are contractions with the second
term. The resulting action $\Ll_\lambda=d_Wi_\lambda+i_\lambda
d_W$ is induced by the co-adjoint action on $\Sg \o\L$,
\cite{GHV}, rel.~(6.5), p.226. With this structure $W(\g)$ becomes
an object of $K(\g)$.

 The cohomology of $W(\g)$ is trivial except in dimension 0,
where it is $k$, \cite{GHV}, Prop.~I, p.228. There are given
canonical maps in $K(\g)$:
\begin{itemize}
\item inclusion $(\Sg )^\g\simeq (\Sg )^\g\o 1\subset W(\g)$,

\item restriction $W(\g)\epi \L$, which sends all the
positive symmetric powers to 0.

\end{itemize}
 The Weil algebra is a model of differential forms on $EG$ and the
sequence of morphisms in $K(\g)$
$$(\Sg)^\g\mono W(g) \epi \L$$
is a model of
$$\Omega^\x(BG)\mono \Omega^\x(EG) \epi \Omega^\x(G)\,.$$

\begin{re}\rm It is easy to see that $W(\g)=\Omega_k\L$. Thus for any
commutative d.g-algebra A $${\rm Hom}_{g-comm}(\L,A)={\rm
Hom}_{d.g-comm}(W(\g),A)\,.$$\end{re}

\section{The distinguished transgression }
 The invariant algebra $(\L)^\g$ is  the exterior algebra
spanned by the space of primitive elements $P^\x$,
 whereas $(\Sg)^\g$ is the symmetric
algebra spanned by some space $\widetilde P^\x$. The point is
that $\widetilde P^\x$ can be canonically chosen.

\begin{pr}\label{trans}{\rm \cite{GHV} Prop.~VI, p.239}. Suppose
$\xi\in P^\x$ is a
primitive element. Then there exist an element $\omega\in W(\g)^\g$
such that
$$ \omega_{|\L}=\xi\,,$$
$$i_\lambda \omega= i_\lambda(1\o\xi)\quad for\, all\;
\lambda\in(\Lg)^\g\,$$
$$d_W(\omega)=\widetilde\xi\o 1\,.$$
\end{pr}

 The element $\omega$ is not unique, but
$\widetilde\xi$ is. The set of $\widetilde\xi$ for $\xi\in P^*$ is
the distinguished space of generators of $(\Sg)^\g$.

\section{Example -- ${\bf su}_2$ }
 The algebra $ {\bf su}_2$ is spanned by ${\bf i}$, ${\bf j}$
and ${\bf k}$ with relation $[{\bf i},{\bf j}]=2{\bf k}$ and its
cyclic transposition. In $\L$ we have $$\dl\i=2\j\wedge\k\quad{\rm
and\;cycl.}$$ In the Weil algebra we have $$d_W(1\o \i)=1\o
2\j\wedge\k +\i\o 1 \quad{\rm and\;cycl.}\,,$$ $$d_W(\i\o
1)=2(\k\o\j-\j\o\k) \quad{\rm and\;cycl.}$$ The primitive elements
in $\L$ are spanned by $\xi=\i\wedge\j\wedge\k$. As $\omega$ of
 Proposition \ref{trans} we take
$$\omega=1\o\i\wedge\j\wedge\k+{\frac12}(\i\o\i+cycl.)$$ Then
$$\widetilde\xi\o 1=d_W(\omega)={\frac12}({\i}^2 +{\j}^2
+{\k}^2)\o 1\,,$$ whereas $$d_W(1\o\xi)=\i\o\j\wedge\k+cycl.$$

\section{Invariant cohomology and equivariant cohomology }

 Denote $(\Lg)^\g$ by $\Lambda_\x$.
Let $D(\Lambda_\x)$ be the derived category of graded
differential $\Lambda_\x$--modules.
For $\M \in K(\g)$ the invariant submodule
$(\M )^\g$ is a differential module over
 $\Lambda_\x$.
We obtain an object in $D(\Lambda_\x)$.
We call it the invariant cohomology of $\M$.

\begin{re}\label{inv} \rm Let $X$ be a manifold on which a compact group
$G$ acts. Let $\g$ be the Lie algebra of $G$. Then $\Omega^\x(X)$ is a
$\g$-module. The invariants of $\Lg$ act on $\Omega^\x(X)$, but this
action does not commute with the differential in general.
We have $[d,i_\lambda]\omega={\cal L}_\lambda\omega$. To obtain
an action which commutes with the differential one
restricts it to $(\Omega^\x(X))^\g$. Fortunately the resulting
cohomology does not
change. We obtain a complex with an action of $\Lambda_\x=H_*(G)$,
which is quasi-isomorphic to $\Omega^\x(X)$. The
cohomology is equal to $H^*(X)$.
\end{re}

 Denote $(\Sg)^\g$ by $S^\x$. Let $D(S^\x)$ be the derived
category of graded differential $S^\x$--modules. Following \cite{Ca}
we define: $$(\M
)_\g:=(\Sg\o \M )^\g$$ with differential $$d_{M,\g}(a\o m)=a\o
d_Mm-\sum_k\lambda^ka\o i_{\lambda_k}m\,.$$ It is a differential
$S^\x$--module. We obtain an object $(\M )_\g$ in $D(S^\x)$. We
call it the equivariant cohomology of $\M $.

 For an object $N^\x$ of $K(\g)$ we define horizontal elements
$$(N^\x)_{hor}=\{n\in N^\x:\;\forall \lambda\in\g\;\;i_\lambda
n=0\}\,.$$ Then define basic elements
$$(N^\x)_{basic}=(N^\x)_{hor}^\g=\{n\in N^\x:\;\forall
\lambda\in\g\;\;i_\lambda n=0,\,i_\lambda dn=0\}\,,$$ which form a
complex.

The following Lemma can be found in \cite{Ca}, but we need to have
an explicit form of the isomorphism, as in \cite{AM} \S4.1.

\begin{lem} \label{Lemma} {\rm\cite{Ca}} The map $\psi_0$ is an
isomorphism of differential graded $S^\x$--modules: $$\psi_0=1\o{\bf T}
_{|1\o\M}:(\M )_\g\rightarrow\left(W(\g)\o \M
\right)_{basic}\,.$$ \end{lem}

{\sc Proof.} The elements of $\left(\L\o \M \right)_{hor}$
are of the form
 $${\bf T}(1\o m)=1\o m -\sum_k\lambda^k\o i_{\lambda_k}m-
\sum_{k<l}\lambda^k\wedge\lambda^l\o
i_{\lambda_k}i_{\lambda_l}m\pm\dots\,,$$ thus they are determined by
$m$. The conclusion follows since $$\left(W(\g)\o \M
\right)_{hor}=\Sg\o \left(\L\o \M \right)_{hor}\,.$$\hfill$\Box$

\begin{re}\rm From the above description we see that
$(\M )_\g$ is an analog of $\Omega^\x(EG\times_GX)$.
Let $X$ and $\g$ be as in \ref{inv}. The construction of the
equivariant cohomology presented here is the so-called  Cartan
model of $\Omega^\x(EG\times_GX)$.
We obtain a complex with an action of $S^\x=H^*(BG)$. The
cohomology is equal to $H^*_G(X)$.
\end{re}

\section{Koszul duality }
 By \cite{GKM}, \S8.5 the following functor
$h:D^{^+}(S^\x)\rightarrow D^{^+}(\Lambda_\x)$ is an equivalence of
categories:
$$h(A^\x)={\rm Hom}_k(\Lambda_\x,A^\x)=(\L)^\g\o A^\x$$
with differential
$$d_h((\xi_1\wedge\dots\wedge\xi_n)\o a)=$$ $$=
\sum_j (-1)^{j+1}(\xi_1\wedge\dots
\vee^j\dots\wedge\xi_n)\o \widetilde\xi_j
a+(-1)^n(\xi_1\wedge\dots\wedge\xi_n)\o da\,,$$
where $\xi_j$'s are primitive.

\begin{th} [Koszul Duality] \label{dual} Let $\g$
be a reductive Lie algebra.
 Suppose $\M$ is an object of $K^{^+}(\g)$
 then in $D^{^+}(\Lambda_\x)$
$$ h((\M )_\g)\simeq (\M )^\g\,.$$
\end{th}

{\sc Proof.} The action of $\g$ on $H^*(W(\g))$ is trivial and,
since $\g$ is reductive, $W(\g)$ is semisimple. Thus by
\cite{GHV}, Th.~V, p.172 the inclusion $$k\o (\M)^\g = W(\g)^\g\o
(\M)^\g \subset (W(\g)\o\M)^\g$$ is a quasi-isomorphism. We want to
construct a quasi-isomorphism $\psi$ from $$h((\M )_\g)=(\L)^\g\o
(\Sg\o \M )^\g $$ to $$(W(\g)\o \M )^\g=(\Sg\o\L\o \M )^\g\,.$$
First we choose a linear map $\omega:P^*\rightarrow W(\g)$
satisfying the conditions of Proposition \ref{trans}. We construct
$\psi$ by the formula extending $\psi_0$ of Lemma \ref{Lemma} with
help of the distinguished transgression of \S5.
$$\psi((\xi_1\wedge\dots\wedge\xi_n)\o m)
=\omega(\xi_1)\dots\omega(\xi_n) \psi_0(m)\,.$$ It is well defined
since $\left(W(\g)\o \M \right)^\g$ is $W(\g)^\g$--module and
$$\psi_0(m)\in\left(W(\g)\o \M \right)_{basic}$$
$$\omega(\xi_1)\in W(\g)^\g\,.$$ The map $\psi$ commutes with
$i_\lambda$ since the image of $\psi_0$ is horizontal. It commutes
with the differential because $d_W(\omega(\xi))=\widetilde\xi\o
1$. This corrects an error in the proof of Lemma 17.6, \cite{GKM},
where $\psi$ does not commute with the differential unless $\g$ is abelian.

We will check that $\psi$ is a quasi-isomorphism. Let's filter both
sides by $S^{ \geq i}(\g^*)$. Then the corresponding quotient
complexes are $$Gr^S_i\psi:(\L)^\g\o (S^i\g^*\o \M
)^\g\longrightarrow (S^i\g^*\o\L\o \M )^\g\,.$$ The differential
on the LHS is just $\epsilon\o 1\o d_M$ (where $\epsilon=(-1)^{{\rm
deg}\xi}$) and the differential on the RHS is $$1\o\dl \o 1
 + \sum_k ad^*_{\lambda_k} \o \lambda^k\wedge\cdot \o 1
+ 1\o \epsilon\o d_M\,.$$ The map $Gr^S_i\psi$ equals
$1\cdot(1\o{\bf T}_{|1\o\M})$. When we untwist it (i.e.~we apply
$\exp(\ii)={\bf T}^{-1}$ to the RHS) the differential takes the form
$$1\o\dl \o 1
 + \sum_k ad^*_{\lambda_k} \o \lambda^k\wedge\cdot \o 1
+ 1\o \epsilon\o d_M+\sum_k 1\o\lambda^k\wedge\cdot \o
\Ll_{\lambda_k}\,.$$ Since we stay in the invariant subcomplex the
differential equals $$1\o\dl \o 1
 - \sum_k 1\o \lambda^k\wedge ad^*_{\lambda_k} \o 1
+ 1\o \epsilon\o d_M\,.$$
 Moreover $\sum_k 1\o \lambda^k\wedge
ad^*_{\lambda_k}=2\dl$, thus the differential on the RHS is
$$-1\o\dl\o 1+1\o \epsilon\o d_M\,.$$ The cohomology of the LHS is
$$H^*\left(\left(\L\right)^\g\o \left(S^i\g^*\o
\M\right)^\g\right) =\left(\L\right)^\g\o H^*\left(\left(S^i\g^*\o
\M \right)^\g\right)$$ and the cohomology of the RHS is $$
H^*\left(\left(S^i\g^*\o \L\o \M \right)^\g\right)=
H^*\left(\left(\L\right)^\g\right)\o H^*\left(\left(S^i\g^*\o
\M\right)^\g\right)$$ again by \cite{GHV}, Th.~V, p.172, since the
action on $H^*(\L)$ is trivial and $\L$ is semisimple. Thus
cohomology of the graded complexes are the same. The conclusion of
\ref{dual} follows. $\Box$

\begin{re}\rm Let $X$ and $\g$ be as in \ref{inv}.
Following \cite{GKM} let us explain the meaning of Koszul duality.
The invariant and the equivariant cohomo\-lo\-gy of $X$ are defined on
the level of derived categories. Theorem \ref{dual} gives a procedure to
reconstruct the invariant cohomology from the equivariant cohomology.
Since this duality is an isomorphism of categories the invariant
cohomology determines the equivariant cohomology as well. The
corresponding statement on the level of graded modules over
$\Lambda_\x$ and $S^\x$ is not true as an easy example in \cite{GKM}
shows. One cannot recover $H^*_G(X)$ from $H^*(X)$
with an action of $H_*(G)$ even in the case $X=S^3$, $G=S^1$. \end{re}

\end{document}